\documentclass[12pt]{article}

\usepackage[centertags]{amsmath}
\usepackage{amsfonts}
\usepackage{amssymb}
\usepackage{amsfonts}
\usepackage{latexsym}

\newcounter{abcd}
\Alph{abcd}

\newtheorem{pred0}{Theorem}

\newcommand {\Hol}{\mathop{\rm Hol}\nolimits}

\newtheorem{prop}{Proposition}
\newtheorem{lemma}{Lemma}

\newtheorem{corol}{Corollary}
\newcommand{\pr}{\noindent{\bf Proof.}\quad }

\newcommand{\epr}{\ $\blacksquare$  \vspace{2mm} }

\newtheorem{rem}{Remark}

\renewcommand{\Im}{\mathop{\rm Im}\nolimits}

\title{Commuting semigroups of holomorphic mappings}

\author{Mark Elin
\\ {\small Department of Mathematics}
\\ {\small ORT  Braude College}
\\ {\small P.O. Box 78, 21982 Karmiel, Israel}
\\ {\small e-mail: mark.elin@gmail.com}
\\ Marina Levenshtein
\\ {\small Department of Mathematics}
\\ {\small The Technion --- Israel Institute of Technology}
\\ {\small 32000 Haifa, Israel}
\\ {\small e-mail: marlev@list.ru}
\\ Simeon Reich
\\ {\small Department of Mathematics}
\\ {\small The Technion --- Israel Institute of Technology}
\\ {\small 32000 Haifa, Israel}
\\ {\small e-mail: sreich@tx.technion.ac.il}
\\ David Shoikhet
\\ {\small Department of Mathematics}
\\ {\small ORT Braude College}
\\ {\small P.O. Box 78, 21982 Karmiel, Israel}
\\ {\small e-mail: davs27@netvision.net.il}}
\date{ }

\begin{document}

\begin{titlepage}

\begin{center}


{\bf\Large Commuting semigroups of holomorphic mappings}

\bigskip

{\large Mark Elin\\}  {\footnotesize Department of Mathematics
\\ ORT  Braude College
\\ P.O. Box 78, 21982 Karmiel, Israel
\\ e-mail: mark.elin@gmail.com}

\bigskip

{\large Marina Levenshtein\\}  {\footnotesize Department of
Mathematics
\\ The Technion --- Israel Institute of Technology
\\ 32000 Haifa, Israel
\\ e-mail: marlev@list.ru}

\bigskip

{\large Simeon Reich\\}  {\footnotesize Department of Mathematics
\\ The Technion --- Israel Institute of Technology
\\ 32000 Haifa, Israel
\\ e-mail: sreich@tx.technion.ac.il}

\bigskip

{\large David Shoikhet\\}  {\footnotesize Department of
Mathematics
\\ ORT  Braude College
\\ P.O. Box 78, 21982 Karmiel, Israel
\\ e-mail: davs27@netvision.net.il}

\end{center}

\bigskip

{\it \noindent 2000 Mathematics Subject Classification: 30C45,
47H10, 47H20. }

\end{titlepage}

\begin{abstract}
Let $S_{1}=\left\{F_t\right\}_{t\geq 0}$ and
$S_{2}=\left\{G_t\right\}_{t\geq 0}$ be two continuous semigroups
of holomorphic self-mappings of the unit disk $\Delta=\{z:|z|<1\}$
generated by $f$ and $g$, respectively. We present conditions on
the behavior of $f$ (or $g$) in a neighborhood of a fixed point of
$S_{1}$ (or $S_{2}$), under which the commutativity of two
elements, say, $F_1$ and $G_1$ of the semigroups implies that the
semigroups commute, i.e., $F_{t}\circ G_{s}=G_{s}\circ F_{t}$ for
all $s,t\geq 0$. As an auxiliary result, we show that the
existence of the (angular or unrestricted) $n$-th derivative of the 
generator $f$ of a semigroup $\left\{F_t\right\}_{t\geq 0}$ 
at a boundary null point of $f$ implies that the corresponding
derivatives of $F_{t}$, $t\geq 0$, also exist, and we obtain
formulae connecting them for $n=2,3$.
\end{abstract}

\section{Introduction}

We denote by $\Hol(\Delta,D)$ the set of all holomorphic functions
on the unit disk $\Delta=\{ z:|z|<1\}$ which map $\Delta$ into a
domain $D\subset\mathbb{C}$, and by $\Hol(\Delta)$ the set of all
holomorphic self-mappings of $\Delta$.

We say that a family $S=\left\{F_t\right\}_{t\geq
0}\subset\Hol(\Delta)$ is {\bf a one-parameter continuous
semigroup on $\Delta$} (a semigroup, in short) if

(i) $F_{t}(F_{s}(z))=F_{t+s}(z)$ for all $t,s\geq 0$

\noindent and

(ii) $\lim\limits_{t\rightarrow 0^+}F_{t}(z)=z$ for all $z\in
\Delta$.

If all the elements $F_{t}$, $t\geq 0$, of a semigroup $S$ are
automorphisms of $\Delta$, then $S$ can be extended to {\bf a
group} of automorphisms $\left\{F_t\right\}_{t\in\mathbb{R}}$ and
property (i) holds for all real $s$ and $t$.

It follows from a result of E. Berkson and H. Porta \cite{B-P}
that each semigroup is differentiable with respect to
$t\in\mathbb{R}^{+}=[0,\infty )$. So, for each one-parameter
continuous semigroup $S=\left\{F_{t}\right\}_{t\geq
0}\subset\Hol(\Delta)$, the limit $$ \lim_{t\rightarrow
0^{+}}\frac{z-F_{t}(z)}{t}=f(z),\quad z\in \Delta, $$  exists and
defines a holomorphic mapping $f\in\Hol(\Delta,\mathbb{C})$. This
mapping $f$ is called the {\bf (infinitesimal) generator of}
$S=\left\{F_{t}\right\}_{t\geq 0}.$ Moreover, the function $
u(t,z):=F_{t}(z), \,  (t,z)\in\mathbb{R}^{+}\times \Delta$,
 is the unique solution of the Cauchy
problem
\begin{equation}\label{2b}
\left\{
\begin{array}{l}
{\displaystyle\frac{\partial u(t,z)}{\partial
t}}+f(u(t,z))=0,\vspace{3mm}\\ u(0,z)=z,\quad z\in \Delta.
\end{array}
\right.
\end{equation}

This solution is univalent on $\Delta$ (see \cite{A}).

We say that $\tau\in\overline{\Delta}$ is a fixed point of
$F\in\Hol(\Delta)$ if either $F(\tau)=\tau$, where $
\tau\in\Delta$, or $\lim\limits_{r\rightarrow 1^-}F(r\tau)=\tau$,
where $\tau\in\partial\Delta=\{z:|z|=1\}$. If $F$ is not an
automorphism of $\Delta$ with an interior fixed point, then by the
Schwarz--Pick Lemma and the Julia--Wolff--Carath\'eodory Theorem,
there is a unique fixed point $\tau\in\overline{\Delta}$ such that
for each $z\in\Delta$,
$\lim\limits_{n\rightarrow\infty}F_{n}(z)=\tau$, where the $n$-th
iteration $F_{n}$ of $F$ is defined by $F_1=F, \, F_{n}=F\circ
F_{n-1},\ n=2,3,\ldots$. Moreover, if $\tau\in\Delta$, then
$|F'(\tau)|<1$, and if $\tau\in\partial\Delta$, then the so-called
angular derivative at the point $\tau$ (see the definition below)
$F'(\tau)\in(0,1]$. This point is called the Denjoy--Wolff point
of $F$. The mapping $F$ is of

--- {\bf dilation type}, if $\tau\in\Delta$,

--- {\bf hyperbolic type}, if $\tau\in\partial\Delta$ and
$0<F'(\tau)<1$,

--- {\bf parabolic type}, if $\tau\in\partial\Delta$ and
$F'(\tau)=1$.

\noindent The mappings of parabolic type fall into two subclasses:

--- {\bf automorphic type}, if all orbits $F_{n}(z)$ are separated
in the hyperbolic Poincar\'{e} metric $\rho$ of $\Delta$, i.e.,
$\lim\limits_{n\rightarrow\infty}\rho (F_{n}(z),F_{n+1}(z))>0 \,
\, \, \mbox{for all} \, \, \, z\in\Delta\,; $

--- {\bf nonautomorphic type}, if no orbit $F_{n}(z)$ is
hyperbolically separated, i.e.,
$\lim\limits_{n\rightarrow\infty}\rho (F_{n}(z),F_{n+1}(z))=0 \,
\, \, \mbox{for all} \, \, \, z\in\Delta\,. $

Consider a semigroup $S=\left\{F_t\right\}_{t\geq 0}$ generated by
$f\in\Hol(\Delta ,\mathbb{C})$. It is a well-known fact that all
elements $F_{t}$ $(t>0)$ of $S$ are of the same type (dilation,
hyperbolic or parabolic) and have the same Denjoy--Wolff point
$\tau$ which is a null point (interior or boundary) of $f$.
(Recall that $\tau\in\partial\Delta$ is a boundary null point of
$f\in\Hol(\Delta , \mathbb{C})$ if $\lim\limits_{r\rightarrow
1^{-}}f(r\tau)=0$.) If $f$ generates a semigroup of dilation type
(which does not consist of automorphisms), then
$\mathrm{Re}f'(\tau)>0$. In the hyperbolic case the angular
derivative $f'(\tau)$ defined by
$f'(\tau):=\lim\limits_{r\rightarrow
1^{-}}\frac{f(r\tau)}{(r-1)\tau}$ exists and is a positive real
number; in the parabolic case $f'(\tau)=0$ (see, for example,
\cite{S}).

We say that a function $f\in\Hol(\Delta,\mathbb{C})$ has an
angular limit $L$ at a point $\tau\in\partial\Delta$ and write
$L:=\angle\lim\limits_{z\rightarrow \tau}f(z)$, if
$f(z)\rightarrow L$ as $z\rightarrow\tau$ in each Stolz angle
$D_{\tau ,\alpha}=\{z\in\Delta :
|\mathrm{arg}(1-\overline{\tau}z)|<\alpha \}, \quad \alpha\in (
0,\frac{\pi}{2})$. If $L$ is finite and the angular limit
$$M:=\angle\lim\limits_{z\rightarrow \tau}\frac{f(z)-L}{z-\tau}$$
exists, then $M$ is said to be the angular derivative $f'(\tau)$.

It is known (see \cite{P}, p. 79) that the existence of the first
angular derivative $f'(\tau)$ of a function $f\in\Hol(\Delta
,\mathbb{C})$ is equivalent to each of the following conditions:

(1) there exists $\angle\lim\limits_{z\rightarrow\tau}f'(z)$, and
then $f'(\tau)=\angle\lim\limits_{z\rightarrow\tau}f'(z)$;

(2) the function $f$ admits the representation
$$f(z)=a_{0}+a_{1}(z-\tau)+\gamma(z),$$ where
$\gamma\in\Hol(\Delta ,\mathbb{C})$,
$\angle\lim\limits_{z\rightarrow \tau}\frac{\gamma(z)}{z-\tau}=0$,
and then $f'(\tau)=a_{1}$.

In Section 2 of this paper we show that higher order angular
derivatives of $f$ can also be defined by either one of these ways
and the definitions are equivalent (Proposition 2). Furthermore,
we show that for a semigroup $\{F_{t}\}_{t\geq 0}$ generated by
$f\in\Hol(\Delta ,\mathbb{C})$, the existence of the $n$-th
($n>1$) angular derivative $f^{(n)}(\tau)$ of $f$ at its boundary
null point $\tau\in\partial\Delta$ implies that for each element
$F_{t}$ of the semigroup, the $n$-th angular derivative at $\tau$
also exists, and obtain formulae connecting $F^{(n)}(\tau)$ with
$f^{n}(\tau)$ for $n=2, 3$ (Theorem 1).

Using these facts, we investigate in Sections 3, 4, and 5
conditions under which the commutativity of two given elements of
the semigroups $S_{1}=\left\{F_t\right\}_{t\geq 0}$ and
$S_{2}=\left\{G_t\right\}_{t\geq 0}$ implies that the semigroups
commute for the dilation, hyperbolic and parabolic cases,
respectively (Theorems 2, 3, and 4).

\section{Higher order boundary derivatives}

We begin by recalling the following known fact.
\begin{prop}[\cite{P}, p. 79]
Let $h$ be holomorphic in $\Delta$. If $\Im h(z)$ has a finite
angular limit at $\tau\in\partial\Delta$, then $(z-\tau)h'(z)$ has
the angular limit $0$ at $\tau$.
\end{prop}

\begin{prop}
Let $f\in\Hol(\Delta,\mathbb{C})$ and let $\tau\in\partial\Delta$.
Then the following assertions are equivalent for any integer
$k\geq 0$:

(i) The function $f$ admits the representation
\begin{equation}\label{d1}
f(z)=\sum_{j=0}^{k}\frac{a_{j}}{j!}(z-\tau) ^{j}+\gamma_{k}(z),
\end{equation}
where $\angle\lim\limits_{z\rightarrow
\tau}\frac{\gamma_{k}(z)}{(z-\tau)^{k}}=0$.

(ii) The angular limit $$\angle\lim\limits_{z\rightarrow
\tau}f^{(k)}(z)$$ exists finitely and coincides with $a_{k}$ in
representation (\ref{d1}).

(iii) For each $0\leq n\leq k$, the angular limit
$$\angle\lim\limits_{z\rightarrow \tau}f^{(n)}(z)$$ exists
finitely and coincides with $a_{n}$ in representation (\ref{d1}).
\end{prop}

\pr

\noindent (i)$\Rightarrow$(ii). Let (i) hold. We show by induction
that for all $0\leq n \leq k$, the following equality is
satisfied:
\begin{equation}\label{d2}
f^{(n)}(z)-\sum_{j=0}^{k-n}\frac{a_{n+j}}{j!}(z-\tau)
^{j}=\gamma_{k-n}(z), \, \, \mbox{where} \, \,
\angle\lim\limits_{z\rightarrow
\tau}\frac{\gamma_{k-n}(z)}{(z-\tau) ^{k-n}}=0.
\end{equation}

It follows from this equality with $n=k$ that
$f^{(k)}(z)-a_{k}=\gamma_{0}(z)$, and, hence,
$\angle\lim\limits_{z\rightarrow \tau}f^{k}(z)=a_{k}$, as
required.

 For $n=0$ relation (\ref{d2}) is obviously equivalent to
(\ref{d1}). Suppose that it holds for $n=m-1\, \, (m\leq k)$,
i.e.,
\begin{equation}\label{d3}
f^{(m-1)}(z)-\sum_{j=0}^{k-m}\frac{a_{m-1+j}}{j!}(z-\tau)
^{j}=\frac{a_{k}}{(k-m+1)!}(z-\tau )^{k-m+1}+\gamma_{k-m+1}(z).
\end{equation}
Denote
$$h(z):=\frac{f^{(m-1)}(z)-\sum_{j=0}^{k-m}\frac{a_{m-1+j}}{j!}(z-\tau)
^{j}}{(z-\tau) ^{k-m+1}}.$$ Then there exists the finite angular
limit
\begin{equation}\label{d4}
\angle\lim\limits_{z\rightarrow \tau} h(z)=\frac{a_{k}}{(k-m+1)!}.
\end{equation}
Now we find
$$(z-\tau)h'(z)=\frac{f^{(m)}(z)-\sum_{j=0}^{k-m-1}\frac{a_{m+j}}{j!}(z-\tau)
^{j} }{(z-\tau) ^{k-m}}-(k-m+1)h(z).$$ Since by Proposition 1,
 $\angle\lim\limits_{z\rightarrow \tau}(z-\tau) h'(z)= 0$, we can
write
$$\frac{f^{(m)}(z)-\sum_{j=0}^{k-m-1}\frac{a_{m+j}}{j!}(z-\tau)
^{j} }{(z-\tau) ^{k-m}}=(k-m+1)h(z)+\mu (z),$$ where
$\angle\lim\limits_{z\rightarrow \tau} \mu (z)=0$. It follows from
this equality and (\ref{d4}), that
$$\frac{f^{(m)}(z)-\sum_{j=0}^{k-m-1}\frac{a_{m+j}}{j!}(z-\tau)
^{j} }{(z-\tau) ^{k-m}}=\frac{a_{k}}{(k-m)!}+\gamma (z),$$ where
$\angle\lim\limits_{z\rightarrow \tau} \gamma (z)=0$. Therefore,
$$f^{(m)}(z)-\sum_{j=0}^{k-m}\frac{a_{m+j}}{j!}(z-\tau) ^{j}
=\gamma_{k-m} (z),$$ where $\gamma_{k-m}(z):=\gamma(z)\cdot
(z-\tau)^{k-m}$ and, consequently,
$\angle\lim\limits_{z\rightarrow
\tau}\frac{\gamma_{k-m}(z)}{(z-\tau)^{k-m}}=0$. In other words,
(\ref{d2}) holds for $n=m$.

\noindent (ii)$\Rightarrow$(iii). Suppose now that there exists
the finite limit
\begin{equation}\label{d41}
a_{k}:=\angle\lim\limits_{z\rightarrow \tau} f^{(k)}(z).
\end{equation}
Consider the equality $$f^{(k-1)}(z)=f^{(k-1)}(0)+\int_{0}^{z}
f^{(k)}(s)ds, \, \, z\in\Delta.$$ Since the angular limit
(\ref{d41}) exists finitely, the function $f^{k}(z)$ is continuous
on each curve $\Gamma (t), \, \, \alpha\leq t \leq\beta , \, \,
\Gamma(\alpha )=0, \, \, \Gamma (\beta )=\tau$, strictly inside
each Stolz angle at $\tau$. Hence, there exists the finite angular
limit $$a_{k-1}:=\angle\lim\limits_{z\rightarrow \tau}
f^{(k-1)}(z)=f^{(k-1}(0)+\int_{0}^{\tau}f^{(k)}(s)ds.$$ Similarly,
for each $0\leq n\leq k$, the limit
\begin{equation}\label{d5}
a_{n}:=\angle\lim\limits_{z\rightarrow \tau} f^{^{(n)}}(z)
\end{equation}
exists finitely.

\noindent (iii)$\Rightarrow$(i). Now we show by induction that for
each $0\leq n\leq k$,
\begin{equation}\label{d6}
f^{(k-n)}(z)=\sum_{j=0}^{n}\frac{a_{k-n+j}}{j!}(z-\tau) ^{j}
+\gamma_{n}(z)
\end{equation}
with $\angle\lim\limits_{z\rightarrow
\tau}\frac{\gamma_{n}(z)}{(z-\tau) ^{n}}=0$.

For $n=0$ equality (\ref{d6}) follows immediately from
(\ref{d41}). Suppose that it holds for $n=m-1 \, \, (m\leq k)$,
i.e.,
\begin{equation}\label{d7}
f^{(k-m+1)}(z)=\sum_{j=0}^{m-1}\frac{a_{k-m+1+j}}{j!}(z-\tau)
^{j}+\gamma_{m-1}(z),
\end{equation}
where $\angle\lim\limits_{z\rightarrow
\tau}\frac{\gamma_{m-1}(z)}{(z-\tau) ^{m-1}}=0$.

It is clear that
$$\frac{f^{(k-m)}(z)-a_{k-m}}{z-\tau}=\int_{0}^{1}f^{(k-m+1)}(t\tau+(1-t)z)dt.$$
Therefore, by (\ref{d5}), $$\angle\lim\limits_{z\rightarrow
\tau}\frac{f^{(k-m)}(z)-a_{k-m}}{z-\tau}=\angle\lim\limits_{z\rightarrow
\tau}\int_{0}^{1}f^{(k-m+1)}(t\tau+(1-t)z)dt=a_{k-m+1}.$$ On the
other hand, by (\ref{d7}),
$$\frac{f^{(k-m)}(z)-a_{k-m}}{z-\tau}=\int_{0}^{1}f^{(k-m+1)}(t\tau+(1-t)z)dt=$$
$$=\int_{0}^{1}\left(\sum_{j=0}^{m-1}\frac{a_{k-m+1+j}}{j!}(t\tau+(1-t)z-\tau)
^{j}+\gamma_{m-1}(t\tau+(1-t)z)\right)dt= $$
$$=\sum_{j=0}^{m-1}\frac{a_{k-m+1+j}}{(j+1)!}(z-\tau)
^{j}+\int_{0}^{1}\gamma_{m-1}(t\tau+(1-t)z)dt. $$
 Hence,
$$f^{(k-m)}(z)=\sum_{j=0}^{m}\frac{a_{k-m+j}}{j!}(z-\tau)
^{j}+\gamma_{m}(z),$$ where
$\gamma_{m}(z)=(z-\tau)\int_{0}^{1}\gamma_{m-1}(t\tau+(1-t)z)dt$.

Now we verify that $\angle\lim\limits_{z\rightarrow
\tau}\frac{\gamma_{m}(z)}{(z-\tau) ^{m}}=0$. Indeed,
$$\angle\lim\limits_{z\rightarrow
\tau}\frac{\gamma_{m}(z)}{(z-\tau)
^{m}}=\angle\lim\limits_{z\rightarrow
\tau}\int_{0}^{1}\frac{\gamma_{m-1}(t\tau+(1-t)z)}{(z-\tau)
^{m-1}}dt=$$ $$=\angle\lim\limits_{z\rightarrow
\tau}\int_{0}^{1}\frac{\gamma_{m-1}(t\tau+(1-t)z)}{(t\tau+(1-t)z-\tau)
^{m-1}} \cdot \frac{(t\tau+(1-t)z-\tau) ^{m-1}}{(z-\tau)
^{m-1}}dt=$$ $$=\int_{0}^{1} \left ((1-t)
^{m-1}\angle\lim\limits_{z\rightarrow
\tau}\frac{\gamma_{m-1}(t\tau+(1-t)z)}{(t\tau +(1-t)z-\tau)
^{m-1}}\right ) dt=0,$$ and for $n=m$ (\ref{d6}) is proved. By
induction, (\ref{d6}) holds for all $0\leq n \leq k$. This
equality with $n=k$ yields representation (\ref{d1}).
\epr

\begin{rem}
It follows from the proof that Proposition 1 also holds if we
replace the angular limit $\angle\lim\limits_{z\rightarrow \tau}$
by the unrestricted limit
$\lim\limits_{\begin{smallmatrix}z\to\tau
\\ z\in\Delta\end{smallmatrix}}$ in (i)--(iii).
\end{rem}

\begin{rem}
Proposition 1 can be rephrased in terms of continuous extension of
the higher order derivatives of $f$ to $\Delta\cup\{\tau\}$
(\cite{B-S}).
\end{rem}

Let $F$ be a holomorphic self-mapping of $\Delta$ and let
$\tau\in\partial\Delta$ be a boundary fixed point of $F$. Then by
the Julia--Wolff--Carath\'eodory Theorem, the first angular
derivative $F'(\tau)$ either exists finitely and is a positive
real number or equals infinity. If $\left\{F_t\right\}_{t\geq 0}$
is a one-parameter continuous semigroup with a boundary fixed
point $\tau\in\partial\Delta$ generated by $f$, then the angular
derivatives $F'_{t}(\tau)$  for all $t>0$ are finite if and only
if the angular derivatives $f'(\tau)=:\beta$ exists finitely.
Moreover, in this case $F'_{t}(\tau)=e^{-\beta t}$ (see
\cite{SD-01}, \cite{E-S-V}, \cite{E-S1}).

As far as the higher order angular derivatives are concerned, even
for the Denjoy--Wolff point one cannot assert that they do exist.
Consider, for example, the parabolic holomorphic self-mapping $F$
of $\Delta$ defined by
\[
F(z):=\frac{2z+(1-z)\mathrm{Log}\left (\frac{2}{1-z} \right )
}{2+(1-z)\mathrm{Log}\left (\frac{2}{1-z} \right )}, \quad
z\in\Delta,
\]
where $\mathrm{Log}$ is the principal branch of the logarithm (see
(\cite{C-M-P})). The Denjoy--Wolff point of this mapping is $\tau
=1$. Consequently, there exists $\angle\lim\limits_{z\rightarrow
1}\frac{\partial F(z)}{\partial z }$. However, the angular limit
$\angle\lim\limits_{z\rightarrow 1}\frac{\partial ^{2}
F(z)}{\partial z ^{2}}$ does not exist finitely.

In Theorem 1 below we show that the existence of the angular
derivatives $f''(\tau)$ and $f'''(\tau)$ of the generator $f$ of a
semigroup $\left\{F_t\right\}_{t\geq 0}$ at a boundary fixed point
$\tau$ implies that for each $t>0$, the angular derivatives
$F''_{t}(\tau):=\angle\lim\limits_{z\rightarrow
\tau}\frac{\partial ^{2} F(z)}{\partial z ^{2}} \,$ and $\,
F'''_{t}(\tau):=\angle\lim\limits_{z\rightarrow
\tau}\frac{\partial ^{3} F(z)}{\partial z ^{3}}$ also exist.
Moreover, we give formulae which connect these derivatives. In the
proof we use the following lemma.

\begin{lemma}[see \cite{POM}, p. 303]
Let $F\in\Hol(\Delta)$ and let $\tau\in\partial\Delta$ be a
boundary fixed point of $F$. If $F$ is conformal at $\tau$, then
nontangential convergence of $z$ to $\tau$ implies that $F(z)$
converges to $\tau$ nontangentially.
\end{lemma}

\begin{pred0}
Let $S=\left\{F_t\right\}_{t\geq 0}$ be a one-parameter continuous
semigroup generated by $f\in\Hol(\Delta , \mathbb{C})$ and let
$\tau\in\partial\Delta$ be a boundary null point of $f$.

(i) If $f'(\tau):=\angle\lim\limits_{z\rightarrow \tau}f'(\tau)$
exists finitely, then for each $t\geq 0$,
$F'_{t}(\tau):=\angle\lim\limits_{z\rightarrow \tau}F'(z)$ also
exists and
\begin{equation}\label{m1}
F'_{t}(\tau) =e^{-\beta t},
\end{equation}
where $\beta=f'(\tau)$.

(ii) If $f''(\tau):=\angle\lim\limits_{z\rightarrow \tau}f''(z)$
exists finitely, then for each $t\geq 0$,
$F''_{t}(\tau):=\angle\lim\limits_{z\rightarrow \tau}F''(z)$ also
exists and
\begin{equation}\label{2y}
F''_{t}(\tau)= \left\{
\begin{array}{l}
-\alpha t, \, \, \, \beta=0 \\ \frac{\alpha}{\beta}e^{-\beta
t}(e^{-\beta t}-1),\, \, \beta\neq 0,
\end{array}
\right.
\end{equation}
where $\beta =f'(\tau)$, $\alpha=f''(\tau)$.

(iii) If $f'''(\tau):=\angle\lim\limits_{z\rightarrow
\tau}f'''(z)$ exists finitely, then for each $t\geq 0$,
$F'''_{t}(\tau):=\angle\lim\limits_{z\rightarrow \tau}F'''(z)$
also exists and
\begin{equation}\label{3y}
F'''_{t}(\tau)= \left\{
\begin{array}{l}
\frac{3}{2}\alpha^{2}t^{2}-\gamma t, \, \, \, \beta=0 \\ \left (
\frac{3\alpha ^{2}}{2\beta ^{2}}+\frac{\gamma}{2\beta}\right )
e^{-3\beta t}-3\frac{\alpha ^{2}}{\beta ^{2}}e^{-2\beta t}+\left (
\frac{3\alpha ^{2}}{2\beta ^{2}}-\frac{\gamma}{2\beta}\right )
e^{-\beta t},\, \, \beta\neq 0,
\end{array}
\right.
\end{equation}
where $\beta =f'(\tau)$, $\alpha=f''(\tau)$, $\gamma=f'''(\tau)$.

\end{pred0}

\pr
Since assertion (i) has been proved in \cite{SD-01} (see also
\cite{C-DM} and \cite{E-S1}), we only present here proofs of
assertions (ii) and (iii).

\noindent (ii) We have already mentioned above that semigroup
elements solve the Cauchy problem (\ref{2b}). Differentiating the
equality
\begin{equation}\label{4x}
\frac{\partial F_{t}(z)}{\partial t}+f(F_{t}(z))=0, \,\,
z\in\Delta, \, \, t\geq 0,
\end{equation}
two times with respect to $z\in\Delta$, we get
\begin{equation}\label{5y}\frac{\partial}{\partial
t}\left(\frac{\partial^{2}F_{t}(z)}{\partial z^{2}}
\right)+f''(F_{t}(z))\left(\frac{\partial F_{t}(z)}{\partial
 z}\right)^{2}+f'(F_{t}(z))\frac{\partial^{2} F_{t}(z)}{\partial
 z^{2}}=0
\end{equation}
for all $z\in\Delta$ and $t\geq 0$.

Define the functions $p(z,t):=f'(F_{t}(z))$,
$q(z,t):=-f''(F_{t}(z)) \left(\frac{\partial F_{t}(z)}{\partial z
} \right) ^{2}$ and
$u_{2}(z,t):=\frac{\partial^{2}F_{t}(z)}{\partial z^{2}}$,
$z\in\Delta$, $t\geq 0$. It is clear that $u_{2}(z,0)=0$.
Rewriting (\ref{5y}) in the form $$\frac{\partial
u_{2}(z,t)}{\partial t}+p(z,t)u_{2}(z,t)=q(z,t), \, \, z\in\Delta,
\, \,  t\geq 0,$$ we find
$$u_{2}(z,t)=e^{-\int_{0}^{t}p(z,s)ds}\cdot
\int_{0}^{t}q(z,s)e^{\int_{0}^{s}p(z,\varsigma)d\varsigma}ds.$$

Now we fix $t$ and let $z$ tend to $\tau$ nontangentially in the
right-hand side of this equality. Since
$\angle\lim\limits_{z\rightarrow\tau}f''(z):=\alpha$ exists
finitely, by Proposition 2, the angular limit
$\angle\lim\limits_{z\rightarrow\tau}f'(z):=\beta$ also exists
finitely. Consequently, for each $t\geq 0$, $\tau$ is a boundary
fixed point of $F_{t}$ and, by item (i),
$\angle\lim\limits_{z\rightarrow\tau}F'_{t}(z)=e^{-\beta t}\neq 0$
(see Theorem 2 in \cite{SD-01}). Hence, by Lemma 1, $F_{t}(z)$
converges to $\tau$ nontangentially as $z$ tends to $\tau$
nontangentially for each $t>0$, and we can conclude that
$\angle\lim\limits_{z\rightarrow\tau}p(z,t)=\beta$ and
$\angle\lim\limits_{z\rightarrow\tau}q(z,t)=-\alpha e^{-2\beta t}$
for each $t>0$. Hence, $$\angle\lim\limits_{z\rightarrow
\tau}\left(e^{-\int_{0}^{t}p(z,t)ds}\cdot\int_{0}^{t}q(z,s)e^{\int_{0}^{s}p(z,\varsigma)d\varsigma}ds
\right)=$$ $$=e^{-\int_{0}^{t}\angle\lim\limits_{z\rightarrow
\tau}p(z,s)ds}\cdot\int_{0}^{t}\angle\lim\limits_{z\rightarrow
\tau}q(z,s)\cdot e^{\int_{0}^{s}\angle\lim\limits_{z\rightarrow
\tau}p(z,\varsigma)d\varsigma}ds=$$ $$=-\alpha e^{-\beta
t}\int_{0}^{t}e^{-\beta s}ds.$$ Therefore if $\beta=0$, then
$$\angle\lim\limits_{z\rightarrow
\tau}\frac{\partial^{2}F_{t}(z)}{\partial z^{2}}=-\alpha t, \, \,
0\leq t<\infty.$$

\noindent If $\beta\neq 0$, then
$$\angle\lim\limits_{z\rightarrow
\tau}\frac{\partial^{2}F_{t}(z)}{\partial
z^{2}}=\frac{\alpha}{\beta}e^{-\beta t}\cdot\left(e^{-\beta
t}-1\right).$$

\noindent (iii) Differentiating equality (\ref{4x}) three times
with respect to $z\in \Delta$, we get $$\frac{\partial}{\partial
t}\left(\frac{\partial^{3}F_{t}(z)}{\partial z^{3}}
\right)+f'''(F_{t}(z))\left(\frac{\partial F_{t}(z)}{\partial
z}\right)^{3}+3f''(F_{t}(z))\frac{\partial F_{t}(z)}{\partial
z}\cdot \frac{\partial^{2} F_{t}(z)}{\partial z^{2}}+$$
\begin{equation}\label{7y} +f'(F_{t}(z))\frac{\partial^{3} F_{t}(z)}{\partial
z^{3}}=0, \, \, t\geq 0,\, \, z\in\Delta.
\end{equation}
Define the functions
$$r(z,t):=-f'''(F_{t}(z))\cdot\left(\frac{\partial
F_{t}(z)}{\partial z}\right)
^{3}-3f''(F_{t}(z))\cdot\frac{\partial F_{t}(z)}{\partial
z}\cdot\frac{\partial^{2}F_{t}(z)}{\partial z^{2}}$$ and
$u_{3}(z,t):=\frac{\partial^{3}F_{t}(z)}{\partial z^{3}}$,
$z\in\Delta$, $t\geq 0$. It is clear that $u_{3}(z,0)=0$.
Rewriting (\ref{7y}) in the form $$\frac{\partial
u_{3}(t,z)}{\partial t}+p(z,t)u_{3}(z,t)=r(z,t), \, \, t\geq 0,$$
we find $$u_{3}(z,t)=e^{-\int_{0}^{t}p(z,s)ds}\cdot
\int_{0}^{t}r(z,s)e^{\int_{0}^{s}p(z,\varsigma)d\varsigma}ds.$$

Now we fix $t$ and let $z$ tend to $\tau$ nontangentially in the
right-hand side of this equality.

Once again, by the continuity of $p(\cdot ,t)$ and $r(\cdot ,t)$ in
$D_{\tau,\nu}\cup \{\tau\}$, $\nu\in(0,\frac{\pi}{2})$,
$$\angle\lim\limits_{z\rightarrow
\tau}\left(e^{-\int_{0}^{t}p(z,s)ds}\cdot\int_{0}^{t}q(z,s)e^{\int_{0}^{s}p(z,\varsigma)d\varsigma}ds\right)=$$
$$=e^{-\int_{0}^{t}\angle\lim\limits_{z\rightarrow
\tau}p(z,s)ds}\cdot\int_{0}^{t}\angle\lim\limits_{z\rightarrow
\tau}q(z,s)\cdot e^{\int_{0}^{s}\angle\lim\limits_{z\rightarrow
\tau}p(z,\varsigma)d\varsigma}ds=$$ $$=-e^{-\beta t
}\cdot\int_{0}^{t}\left(\gamma e^{-3\beta s}+3\alpha e^{-\beta
s}\cdot \angle\lim\limits_{z\rightarrow
\tau}\frac{\partial^{2}F_{s}(z)}{\partial z^{2}} \right) e^{\beta
s }ds.$$

By Proposition 2, the limit $\angle\lim\limits_{z\rightarrow
\tau}\frac{\partial^{2}F_{t}(z)}{\partial z^{2}}$ exists and by
item (ii) proved above, it is given by equality (\ref{2y}).

Hence, the limit $\angle\lim\limits_{z\rightarrow
\tau}\frac{\partial^{3}F_{t}(z)}{\partial z^{3}}$ exists and in
the parabolic case $(\beta=0)$ it equals
$$\angle\lim\limits_{z\rightarrow
\tau}\frac{\partial^{3}F_{t}(z)}{\partial
z^{3}}=-\int_{0}^{t}(\gamma-3\alpha^{2}s)ds=\frac{3\alpha^{2}t^{2}}{2}-\gamma
t.$$ In the hyperbolic case $(\beta\neq 0)$ this limit also exists
and  $$\angle\lim\limits_{z\rightarrow
\tau}\frac{\partial^{3}F_{t}(z)}{\partial z^{3}}=-e^{-\beta
t}\cdot\int_{0}^{t}\left( \left( \gamma +\frac{3\alpha^{2}}{\beta}
\right)e^{-2\beta s}-\frac{3\alpha^{2}}{\beta}e^{-\beta s}\right)
ds=$$  $$=\left ( \frac{3\alpha ^{2}}{2\beta
^{2}}+\frac{\gamma}{2\beta}\right ) e^{-3\beta t}-3\frac{\alpha
^{2}}{\beta ^{2}}e^{-2\beta t}+\left ( \frac{3\alpha ^{2}}{2\beta
^{2}}-\frac{\gamma}{2\beta}\right ) e^{-\beta t}.$$ \epr

\begin{corol}
Let $f\in \Hol(\Delta,\mathbb{C})$ be the generator of a parabolic
semigroup $\left\{F_t\right\}_{t\geq 0}$ with the Denjoy--Wolff
point $\tau\in\partial\Delta$. If $\angle\lim\limits_{z\rightarrow
\tau}f''(\tau)=\angle\lim\limits_{z\rightarrow \tau}f'''(\tau)=0$,
then $F_{t}=\mathrm{I}$ for all $t\geq 0$.
\end{corol}

Indeed, these conditions imply that $F'_{t}(\tau)=1, \, \,
F''_{t}(\tau)=F'''_{t}(\tau)=0$ for all $t\geq 0$ and, by
\cite{E-L-S-T}, we get $F_{t}=\mathrm{I}$.

\begin{rem} \emph{As a matter of fact, repeating our proof and using Remark 1, one can show that the
angular limits in Theorem 1 can be replaced by unrestricted
limits. Namely:}

Let $S=\left\{F_t\right\}_{t\geq 0}$ be the semigroup generated by
$f$. Assume that for each $t>0$ the unrestricted limit $\lim\limits_{\begin{smallmatrix}z\to \tau \\
z\in\Delta\end{smallmatrix}}F(z)$ exists, where $\tau$ is a
boundary null point of $f$. The following assertions hold:

(i) If the unrestricted limit
$\beta:=\lim\limits_{\begin{smallmatrix}z\to \tau \\
z\in\Delta\end{smallmatrix}}f'(z)$ exists finitely, then
$$\lim\limits_{\begin{smallmatrix}z\to \tau \\
z\in\Delta\end{smallmatrix}}F'_{t}(z)=e^{-\beta t} \, \, \,   for
\, \,  each \, \,  t\geq 0. $$

(ii) If the unrestricted limit
$\alpha:=\lim\limits_{\begin{smallmatrix}z\to \tau \\
z\in\Delta\end{smallmatrix}}f''(z)$ exists finitely, then
\begin{equation}\label{20y}
\lim\limits_{\begin{smallmatrix}z\to \tau \\
z\in\Delta\end{smallmatrix}}F''_{t}(z)= \left\{
\begin{array}{l}
-\alpha t, \, \, \, \beta=0 \\ \frac{\alpha}{\beta}e^{-\beta
t}(e^{-\beta t}-1),\, \, \beta\neq 0,
\end{array}
\right.
\end{equation}
for each $t\geq 0$.

(iii) If the unrestricted limit
$\gamma:=\lim\limits_{\begin{smallmatrix}z\to \tau \\
z\in\Delta\end{smallmatrix}}f'''(z)$ exists finitely, then
\begin{equation}\label{r30y}
\lim\limits_{\begin{smallmatrix}z\to \tau \\
z\in\Delta\end{smallmatrix}}F'''_{t}(z)= \left\{
\begin{array}{l}
\frac{3}{2}\alpha^{2}t^{2}-\gamma t, \, \, \, \beta=0 \\ \left (
\frac{3\alpha ^{2}}{2\beta ^{2}}+\frac{\gamma}{2\beta}\right )
e^{-3\beta t}-3\frac{\alpha ^{2}}{\beta ^{2}}e^{-2\beta t}+\left (
\frac{3\alpha ^{2}}{2\beta ^{2}}-\frac{\gamma}{2\beta}\right )
e^{-\beta t},\, \, \beta\neq 0,
\end{array}
\right.
\end{equation}
for each $t\geq 0$.
\end{rem}

\begin{rem}
The arguments used in the proof of Theorem 1 can be used to derive
analogous results for derivatives of any order $k\geq 4$.

\end{rem}

\section{Semigroups with an interior fixed point}

In our proofs we use the two following facts established by C. C.
Cowen in \cite{C}.

\begin{prop}
Let $F$, $G_{1}$, $G_{2}$ be holomorphic self-mappings of
$\Delta$, not automorphisms of $\Delta$, and let $G_{1}$ and
$G_{2}$ commute with $F$. Suppose that $\tau\in\overline{\Delta}$
is the Denjoy--Wolff point of $F$ and that $0<|F'(\tau)| <1$. Then
$G_{1}$ and $G_{2}$ commute with each other.
\end{prop}

\begin{prop}
Let $F$ and $G$ be two commuting holomorphic self-mappings of
$\Delta$, not automorphisms of $\Delta$, and let
$\tau\in\overline{\Delta}$ be their common Denjoy--Wolff point.

(i) If $F'(\tau)=0$, then $G'(\tau)=0$.

(ii) If  $0<|F'(\tau)|<1$, then $0<|G'(\tau)|<1$.

(iii) If $F'(\tau)=1$, then $G'(\tau)=1$.
\end{prop}

The following fact is more or less known (see, for example,
\cite{A}).
\begin{prop}
Let $S=\left\{F_{t}\right\}_{t\geq 0}$ be a semigroup in $\Delta$.
Assume $F_{t_{0}}$ is an automorphism of $\Delta$ for some
$t_{0}>0$; then each element $F_{t}$ of $S$ is an automorphism of
$\Delta$.
\end{prop}

We now begin our investigation of commuting semigroups. Note that
in all the following theorems the condition $F_{1}\circ
G_{1}=G_{1}\circ F_{1}$ can be replaced by the condition
$F_{p}\circ G_{q}=G_{q}\circ F_{p}$ for some $p,q>0$.

\begin{pred0}[dilation case]
Let $S_{1}=\left\{F_{t}\right\}_{t\geq 0}$ and
$S_{2}=\left\{G_{t}\right\}_{t\geq 0}$ be two continuous
semigroups on $\Delta$ generated by $f$ and $g$, respectively, and
let $F_{1}\circ G_{1}=G_{1}\circ F_{1}$. Suppose that $f$ has an
interior null point $\tau\in\Delta$.

(i) If $S_{1}$ and $S_{2}$ are not groups of automorphisms of
$\Delta$, then they commute.

(ii) If $S_{1}$ is a nontrivial group of elliptic automorphisms of
$\Delta$ and $S_{2}$ is a semigroup of self-mappings of $\Delta$,
then $S_{1}$ and $S_{2}$ commute if and only if  $S_{2}$ is a
semigroup of linear fractional transformations of the form
\begin{equation}\label{201b}
G_{t}(z)=m_{\tau}(e^{-at}\cdot m_{\tau}(z))
\end{equation}
for some $a\in\mathbb{C}$, where
$m_{\tau}(z)=\frac{\tau-z}{1-\overline{\tau}z}$.
\end{pred0}

Note that the function $G_{t}$ defined by equality (\ref{201b}) is
a self-mapping of $\Delta$ if and only if $\mathrm{Re}\, a\geq 0$.

\pr Since $\tau$ is an interior null point of the generator $f$,
it is the unique interior fixed point of the semigroup $S_{1}$
(see \cite{A}). The commutativity of $F_{1}$ and $G_{1}$ implies
that $\tau$ is a fixed point of $G_{1}$ and, consequently, $\tau$
is a fixed point of $G_{t}$ for each $t>0$.

(i) If $S_{1}$ and $S_{2}$ are not groups of automorphisms of
$\Delta$, then $0<|F'_{t}(\tau)|<1$ and $0<|G'_{t}(\tau)|<1$ for
all $t>0$, by the Schwarz--Pick Lemma and the univalence of
$F_{t}$ and $G_{t}$ on $\Delta$ for all $t\geq 0$.

The function $G_{1}$ commutes with $F_{1}$ (by our assumption)
and, for each $t\geq 0$, the mapping $F_{t}$ commutes with $F_{1}$
(by the semigroups property). Therefore Proposition 3 implies
that $G_{1}\circ F_{t}=F_{t}\circ G_{1}$ for all $t\geq 0$.

Fix an arbitrary $t>0$. Similarly, since $G_{1}\circ
F_{t}=F_{t}\circ G_{1}$ and $G_{1}\circ G_{s}=G_{s}\circ G_{1}$
for all $s\geq 0$, we get, by Proposition 3, that $G_{s}\circ
F_{t}=F_{t}\circ G_{s}$ for all $s\geq 0$. Hence, the semigroups
$S_{1}$ and $S_{2}$ commute, as claimed.

(ii) Since $S_{1}$ is a group of elliptic automorphisms of
$\Delta$ with a fixed point $\tau\in\Delta$, the functions $F_{t}$
are of the form (see \cite{B-K-P}) $$F_{t}(z)=m_{\tau}(e^{i\varphi
t} m_{\tau}(z)) \, \, \, \mbox{for some } \, \,
\varphi\in\mathbb{R}.$$

Let $G_{t}(z)=m_{\tau}(e^{-at} m_{\tau}(z))$. Using the equality
$m_{\tau}(m_{\tau}(z))=z$, we get
$$F_{t}(G_{s}(z))=m_{\tau}(e^{i\varphi t}
m_{\tau}(m_{\tau}(e^{-as} m_{\tau}(z))) )=m_{\tau}(e^{i\varphi t}
e^{-as} m_{\tau}(z))=$$ $$=m_{\tau}(e^{-as} e^{i\varphi t}
m_{\tau}(z))=m_{\tau}(e^{-as} m_{\tau}(m_{\tau}(e^{i\varphi t}
m_{\tau}(z) )))=G_{s}(F_{t}(z)).$$

Conversely, suppose that $F_{t}\circ G_{s}=G_{s}\circ F_{t}$ for
all $s,t\geq 0$. Denote $$\widetilde{F}_{t}(z)=e^{i\varphi t}z, \,
\widetilde{G}_{t}=m_{\tau}\circ G_{t}\circ m_{\tau}.$$ Then
$\{\widetilde{F}_{t}\}_{t\geq 0}$ is a group of automorphisms of
$\Delta$ with a fixed point at zero, and
$\{\widetilde{G}_{t}\}_{t\geq 0}$ is a semigroup of self-mappings
of $\Delta$ with a fixed point at zero. It is obvious that the
semigroups $\{\widetilde{F}_{t}\}_{t\geq 0}$ and
$\{\widetilde{G}_{t}\}_{t\geq 0}$ commute. Consequently, their
generators $\widetilde{g}(z)$ and $\widetilde{f}(z)=-i\varphi z$
are proportional (see \cite{E-L-S-T}). So $\widetilde{g}(z)=az$
for some $a\in\mathbb{C}$. Therefore $\widetilde{G}(z)=e^{-at} z$
and $G_{t}(z)=m_{\tau}(e^{-at} m_{\tau}(z))$. \epr

We see from this theorem that if $S_{1}$ is a group of elliptic
automorphisms, the commutativity of $F_{1}$ and $G_{1}$ does not
imply that the semigroups $S_{1}$ and $S_{2}$ commute.
Nevertheless, in this case one can still obtain some additional
information about the semigroup $S_{2}$. The following assertions
explain our claim.

\begin{prop}
If $S_{1}=\left\{F_{t}\right\}_{t\geq 0}$ is a group of elliptic
automorphims whereas $S_{2}=\left\{G_{t}\right\}_{t\geq 0}$ is a
semigroup of self-mappings of $\Delta$ which are not
automorphisms, then the commutativity of $F_{1}$ and $G_{1}$
implies that $F_{1}\circ G_{t}=G_{t}\circ F_{1}$ for all $t\geq
0$.
\end{prop}

\pr
Let $\tau\in\Delta$ be the common fixed point of $S_{1}$ and
$S_{2}$. Then the functions $F_{t}$ are of the form
$F_{t}(z)=m_{\tau}(e^{i\varphi t} m_{\tau}(z)), \, \,
\varphi\in\mathbb{R}, \, \, z\in\Delta$, where
$m_{\tau}(z)=\frac{\tau-z}{1-\overline{\tau}z}$.

Denote $\widetilde{F}_{t}(z)=e^{i\varphi t} z$ and
$\widetilde{G}_{t}(z)=m_{\tau}(G_{t}(m_{\tau}(z)))$. Then
$\{\widetilde{F}_{t}\}_{t\geq 0}$ is a group of automorphisms of
$\Delta$ with its common fixed point at zero, and
$\{\widetilde{G}_{t}\}_{t\geq 0}$ is a semigroup of self-mappings
of $\Delta$ which are not automorphisms with its common fixed
point also at zero.

It is obvious that for each $t>0$, $F_{1}$ and $G_{t}$ commute if
and only if $\widetilde{F}_{1}$ and $\widetilde{G}_{t}$ commute.
Hence, by our assumption, $\widetilde{F}_{1}\circ
\widetilde{G}_{1}=\widetilde{G}_{1}\circ \widetilde{F}_{1}$ or,
which is the same,
$e^{i\varphi}\widetilde{G}_{1}(z)=\widetilde{G}_{1}(e^{i\varphi}z)$.
It follows that for all $n\in\mathbb{N}$, $\widetilde{F}_{1}\circ
\widetilde{G}_{n}=\widetilde{G}_{n}\circ \widetilde{F}_{1}$, where
$\widetilde{G}_{n}$ are the iterates of $\widetilde{G}_{1}$, i.e.,
$ \widetilde{G}_{n}=\widetilde{G}_{1}\circ \widetilde{G}_{n-1}$.

Since $\widetilde{G}_{1}$ is a self-mapping of $\Delta$ (which is
not an automorphism) with a fixed point at the origin, there
exists a unique univalent solution $h$ of the functional equation
$$h(\widetilde{G}_{1}(z))=\alpha h(z), \, \, \mbox{with} \, \,
\alpha=\widetilde{G}'_{1}(0),$$ normalized by $h(0)=0, \, \,
h'(0)=1$ (see, for example, \cite{Shap}). This solution is given
by
$$h(z)=\lim\limits_{n\rightarrow\infty}\frac{\widetilde{G}_{n}(z)}{\alpha^{n}}.$$
Moreover, by \cite{E-G-R-S}, for all real positive $t$,
$$h(\widetilde{G}_{t}(z))=\alpha^{t}h(z) \, .$$ Therefore,
$$h(\widetilde{F}_{1}(\widetilde{G}_{t}(z)))=h(e^{i\varphi}\widetilde{G}_{t}(z))=\lim\limits_{n\rightarrow\infty}\frac{\widetilde{G}_{n}(e^{i\varphi}\widetilde{G}_{t}(z))}{\alpha^{n}}=\lim\limits_{n\rightarrow\infty}\frac{e^{i\varphi}\widetilde{G}_{n}(\widetilde{G}_{t}(z))}{\alpha^{n}}=$$
$$=e^{i\varphi}h(\widetilde{G}_{t}(z))=e^{i\varphi}\alpha^{t}h(z)=\alpha^{t}\lim\limits_{n\rightarrow\infty}\frac{e^{i\varphi}\widetilde{G}_{n}(z)}{\alpha^{n}}=\alpha^{t}\lim\limits_{n\rightarrow\infty}\frac{\widetilde{G}_{n}(e^{i\varphi}z)}{\alpha^{n}}=$$
$$=\alpha^{t}h(e^{i\varphi}z)=h(\widetilde{G}_{t}(e^{i\varphi}z))=h(\widetilde{G}_{t}(\widetilde{F}_{1}(z)))$$
and, by the univalence of $h$, we get $\widetilde{F}_{1}\circ
\widetilde{G}_{t}=\widetilde{G}_{t}\circ \widetilde{F}_{1}$ for
all $t\geq 0$. Consequently, $F_{1}$ and $G_{t}$ commute for all
$t\geq 0$. \epr

\begin{corol}
Let $S_{1}=\left\{F_{t}\right\}_{t\geq 0}$ be a group of elliptic
automorphisms of $\Delta$, i.e., $F_{t}(z)=m_{\tau}(e^{i\varphi
t}m_{\tau}(z))$, $\varphi\in\mathbb{R}$, $\tau\in\Delta$, and let
$S_{2}=\{G_{t} \}_{t\geq 0}$ be a semigroup of self-mappings of
$\Delta$. Suppose that $\frac{\varphi}{\pi}$ is an irrational
number and $F_{1}$ and $G_{1}$ commute. Then
$G_{t}(z)=m_{\tau}(e^{-at}m_{\tau}(z)), \, \, a\in\mathbb{C},$
and, consequently, the semigroups $S_{1}$ and $S_{2}$ commute.
\end{corol}

\pr
Once again, we define the functions $\widetilde{F}_{t}=e^{i\varphi
t}z$ and $\widetilde{G}_{t}=m_{\tau}\circ G_{t}\circ m_{\tau}$.
The commutativity of $F_{1}$ and $G_{1}$ implies that
$\widetilde{F}_{1}\circ \widetilde{G}_{1}=\widetilde{G}_{1}\circ
\widetilde{F}_{1}$ and, by Proposition 6, $\widetilde{F}_{1}\circ
\widetilde{G}_{t}=\widetilde{G}_{t}\circ \widetilde{F}_{1}$ for
all $t\geq 0 $. Therefore
$\widetilde{G}_{t}(e^{in\varphi}z)=e^{in\varphi}\widetilde{G}_{t}(z)$
for all $n\in\mathbb{\mathbb{N}}$. Since the set $\{e^{in\varphi}
\}_{n\in\mathbb{N}}$ is dense in the unit circle,
$\widetilde{G}_{t}(\lambda z)=\lambda \widetilde{G}_{t}(z)$ for
all $\lambda$ with $|\lambda|=1$ and $z\in\Delta$, by the
continuity of $\widetilde{G}_{t}$ on $\Delta$.

Fix $0\neq z\in\Delta$ and $t>0$, and consider the analytic
function $q(\lambda)$ on the closed unit disk defined by

\begin{equation}\label{200b}
q(\lambda)=\left\{
\begin{array}{l}
{\displaystyle\frac{\widetilde{G}_{t}(\lambda z)}{\lambda}}, \, \,
\lambda\neq 0 ,\vspace{3mm}\\
{\displaystyle\lim\limits_{\lambda\rightarrow
0}\frac{\widetilde{G}_{t}(\lambda z)}{\lambda}=\left .
z\frac{\partial}{\partial w
}\widetilde{G}_{t}(w)\right|_{w=0}},\quad \lambda=0.
\end{array}
\right.
\end{equation}

This function is constant on the unit circle:
$q(\lambda)=\widetilde{G}_{t}(z)$. Moreover, $q(\lambda)\neq 0$
for all $\lambda\in\Delta$. Therefore
$q(\lambda)=\widetilde{G}_{t}(z)$ for all
$\lambda\in\overline{\Delta}$. So for each $z\neq 0$ and $t>0$,
$\widetilde{G}_{t}(\lambda z)=\lambda \widetilde{G}_{t}(z) $.
Consequently, this equality holds for all $z\in\Delta$. Hence
$\widetilde{G}_{t}$ is a linear function for each $t>0$, i.e.,
$\widetilde{G}_{t}(z)=e^{-at}z$ for some $a\in\mathbb{C}$,
$\mathrm{Re} \, a\geq 0$, and the assertion follows. \epr

In contrast with this corollary, if $\frac{\varphi}{\pi}$ is a
rational number, the semigroups $S_{1}$ and $S_{2}$ do not
necessarily commute. The following example gives a large class of
semigroups $S_{2}=\{G_{t} \}_{t\geq 0}$ such that $F_{1}\circ
G_{t}=G_{t}\circ F_{1}$ for all $t\geq 0$, but the semigroups
$S_{1}$ and $S_{2}$ do not commute.

{\vspace{2mm}\bf Example.} Let $S_{1}=\left\{F_{t}\right\}_{t\geq
0}$, where $F_{t}(z)=e^{i\frac{2\pi}{n}t}z, \, n\in\mathbb{N},$
and let $S_{2}=\left\{G_{t}\right\}_{t\geq 0}$ be the semigroup
generated by $g(z)=z p(z^{n})$, where $\mathrm{Re} \, p(z^{n})\geq
0$ for all $z\in\Delta$. Then $F_{1}\circ G_{t}=G_{t}\circ F_{1}$
for all $t\geq 0$.

Indeed, denote $u=u(t,z):=G_{t}(z)$. Then $u$ is the unique
 solution of the Cauchy problem
\begin{equation}\label{z01}
\left\{
\begin{array}{l}
{\displaystyle\frac{\partial u}{\partial
t}}+up(u^{n})=0,\vspace{3mm}\\ u(0,z)=z,\quad z\in \Delta,
\end{array}
\right.
\end{equation}
and, consequently,
\begin{equation}\label{z02}
\int_{z}^{G_{t}(z)}\frac{d\varsigma}{\varsigma p(\varsigma
^{n})}=-t \, \, \, \,  \mbox{for all} \, \, \, \,  z\in\Delta.
\end{equation}
Substituting $e^{i\frac{2\pi}{n}}z$ instead of $z$, we get
\[
\int_{e^{i\frac{2\pi}{n}}z}^{G_{t}(e^{i\frac{2\pi}{n}}z)}\frac{d
\varsigma}{\varsigma p(\varsigma ^{n})}=-t.
\]
Now substitute $\varsigma =e^{i\frac{2\pi}{n}}w$:
\begin{equation}\label{z03}
\int_{z}^{G_{t}(e^{i\frac{2\pi}{n}}z)e^{-i\frac{2\pi}{n}}}\frac{d
w}{w p(w^{n}e^{i 2 \pi
})}=\int_{z}^{G_{t}(e^{i\frac{2\pi}{n}}z)e^{-i\frac{2\pi}{n}}}\frac{d
w}{w p(w^{n})}=-t, \, \,  z\in\Delta.
\end{equation}

Equalities (\ref{z02}) and (\ref{z03}) imply that
\begin{equation}\label{z201}
\int_{G_{t}(z)}^{G_{t}(e^{i\frac{2\pi}{n}}z)e^{-i\frac{2\pi}{n}}}\frac{d
w}{w p(w^{n})}=0, \, \,  z\in\Delta.
\end{equation}

By the uniqueness of the solution to the Cauchy problem
(\ref{z01}), the equation $$\int_{z}^{u}\frac{d w}{w p(w^{n})}=-s,
\, \, s\geq 0, \, \, z\in\Delta,$$ has the unique solution
$u=G_{s}(z)$ for each $s\geq 0$. Thus, it follows from
(\ref{z201}) that
$G_{t}(e^{i\frac{2\pi}{n}}z)e^{-i\frac{2\pi}{n}}=G_{0}(G_{t}(z))=G_{t}(z)$.
Hence, $G_{t}(e^{i\frac{2\pi}{n}}z)=e^{i\frac{2\pi}{n}}G_{t}(z)$.
Therefore $F_{1}$ commutes with $G_{t}$ for all $t\geq 0$. At the
same time, if $p$ is not a constant function, the semigroups do
not commute because their generators are not proportional.

\section{Semigroups of hyperbolic type}

We start this section with an assertion which is of independent
interest.

\begin{prop}
Let $F$ and $G$ be two commuting holomorphic self-mappings of
$\Delta$ and assume that $G$ is not the identity. If $F$ is of hyperbolic 
type, then $G$ is of hyperbolic type too.
\end{prop}

\pr If $F$ is a hyperbolic automorphism of $\Delta$, then by
Lemma~2.1 in \cite{H} $G$ is a hyperbolic automorphism of
$\Delta$.

Let $F$ be a holomorphic self-mapping of $\Delta$ which is not an
automorphism of $\Delta$. In this case, by a result in \cite{Beh},
the mappings $F$ and $G$ have a common Denjoy--Wolff point
$\tau\in\partial\Delta$. We have to show that $G$ is of hyperbolic
type, i.e., $0<G'(\tau)<1$. Suppose, to the contrary, that $G$ is
of parabolic type, i.e., $G'(\tau)=1$. Then, by Proposition 4(ii),
$G$ must be a parabolic automorphism.

Denote $g:=C\circ G\circ C^{-1}$ and $f:=C\circ F\circ C^{-1}$,
where $C(z)=\frac{\tau +z}{\tau -z}$. Then $f$ and $g$ are two
commuting holomorphic self-mappings of the right half-plane
$\mathbb{H}=\{z:\mathrm{Re}z>0\}$ with their common Denjoy--Wolff
point at infinity. Moreover, $g$ is a parabolic automorphism of
$\mathbb{H}$ while $f$ is a hyperbolic self-mapping of
$\mathbb{H}$. Consequently, $f$ and $g$ are of the forms (see
\cite{Shap}): $$f(w)=cw+\Gamma_{F}(w) \quad \mbox{with} \quad
c=\frac{1}{F'(\tau)}>1 \quad \mbox{and} \quad
\angle\lim\limits_{w\rightarrow\infty}\frac{\Gamma _{F}(w)}{w}=0,$$
and
$$g(w)=w+ib \quad\mbox{with} \quad b\in\mathbb{R}\setminus \{0\}
\quad \mbox{and} \quad w\in\mathbb{H}.$$ By a simple
calculation and the commutativity of $f$ and $g_{n}$, we infer
from the above representations that
\begin{equation}\label{z24}
f(w+nib)=f(w)+nib, \quad w\in\mathbb{H}.
\end{equation}
Hence,
\[
\frac{f(w+nib)}{w+nib}=\frac{f(w)}{w+nib}+\frac{nib}{w+nib}\,,
\quad w\in\mathbb{H}.
\]
Letting $n\rightarrow\infty$, we obtain that for each
$w\in\mathbb{H}$, the limit
$\lim\limits_{n\rightarrow\infty}\frac{f(w+nib)}{w+nib}$ exists
and equals 1.

Fix $w_{0}\in\mathbb{H}$. Consider the curve $l:=\{w_{0}+it:\,
t\in\mathbb{R}, \ \mathop{\rm sgn}t=\mathop{\rm sgn}b\}$. We
intend to show that the limit $\lim\limits_{l\ni
z\rightarrow\infty}\frac{f(z)}{z}$ exists and equals 1.

To this end, fix an arbitrary $\varepsilon>0$ and take
$N\in\mathbb{N}$ such that
\[
N>\frac{1}{|b|}\left(\frac{|f(w)-w|}{\varepsilon}+|w|\right)\
\mbox{ and }\ N>\frac{|w|}{|b|}
\]
for all $w\in [w_{0},w_{0}+ib]$.

Then $\left|\frac{f(z)}{z}-1\right|< \varepsilon$ for all $z\in l$
with $\mathop{\rm sgn}b\cdot \Im z>\mathop{\rm sgn}b(\Im
w_{0}+Nb)$.

Indeed, if $\mathop{\rm sgn}b\cdot \Im z>\mathop{\rm sgn}b(\Im
w_{0}+Nb)$, then $z=\alpha +ikb$ for some $\alpha\in
[w_{0},w_{0}+ib]$ and $k\geq N$.

Hence, $k|b|\geq |\alpha|$ and
$k>\frac{1}{|b|}\left(\frac{|f(\alpha)-\alpha|}{\varepsilon}+|\alpha|\right)$.
Consequently, $|\alpha
+ikb|>k|b|-|\alpha|>\frac{|f(\alpha)-\alpha|}{\varepsilon}$.

Now using (\ref{z24}), we obtain that
$$\left|\frac{f(z)}{z}-1\right|=\left|\frac{f(\alpha +kib)
}{\alpha +kib}-1\right|=\left|\frac{f(\alpha)-\alpha}{\alpha +ikb
}\right|<\varepsilon .$$

Thus $\lim\limits_{l\ni
z\rightarrow\infty}\frac{f(z)}{z}=1$. It now follows from  
Lindel\"{o}f's theorem (see, for example, \cite{S}) that
$\angle\lim\limits_{z\rightarrow\infty}\frac{f(z)}{z}=1$, which
contradicts our assumption. Therefore the mapping $G$ is indeed of
hyperbolic type. \epr

\begin{pred0}[hyperbolic case]
Let $S_{1}=\left\{F_{t}\right\}_{t\geq 0}$ and
$S_{2}=\left\{G_{t}\right\}_{t\geq 0}$ be continuous semigroups on
$\Delta$ generated by $f$ and $g$, respectively, and assume that
$F_{1}\circ G_{1}=G_{1}\circ F_{1}$. Suppose that $f$ has a
boundary null point $\tau\in\partial\Delta$, such that
$f'(\tau):=\angle\lim\limits_{z\rightarrow\tau}f'(z)>0$, i.e., the
semigroup $S_{1}$ is of hyperbolic type. Then the semigroups
$S_{1}$ and $S_{2}$ commute. Thus, if $g\neq 0$ then $S_{2}$ is
also of hyperbolic type.
\end{pred0}

\pr
By our assumption, $\tau$ is the Denjoy--Wolff point of the
semigroup $S_{1}$.

First we suppose that $S_{1}$ and $S_{2}$ consist of automorphisms
of $\Delta$. Since $f'(\tau)>0$, $S_{1}$ consists of hyperbolic
automorphisms of $\Delta$ and its generator $f$ is of the form
$$f(z)=\frac{a_{1}}{\tau-\varsigma}(z-\tau)(z-\varsigma),$$ where
$a_{1}$ is a positive real number and $\varsigma$ is the second
common fixed point of the semigroup $S_{1}$ (see \cite{B-K-P}).

The commutativity of $F_{1}$ and $G_{1}$ implies that $G_{1}$ has
the same fixed points $\tau$ and $\varsigma$; consequently,
$S_{2}$ consists of hyperbolic automorphisms of $\Delta$, and its
generator $g$ is of the form
$$g(z)=\frac{a_{2}}{\varsigma-\tau}(z-\tau)(z-\varsigma),$$ where
$a_{2}$ is a non-zero real number. Hence, $g(z)= -
\frac{a_{2}}{a_{1}}f(z)$, and by Theorem 3 in \cite{E-L-S-T}, the
semigroups commute.

Suppose now that at least one of the semigroups $S_{1}$ and
$S_{2}$ consists of self-mappings of $\Delta$ which are not
automorphisms. By a result in \cite{Beh}, $\tau$ is the common
Denjoy--Wolff point of $S_{1}$ and $S_{2}$. Moreover, by
Theorem~1, $\alpha:=F'_{1}(\tau)=e^{-f'(\tau)}\in(0,1)$.
Consequently, by Proposition~7, $\beta:=G'_{1}(\tau)\in(0,1)$.

Since $F_{1}$ is a hyperbolic self-mapping of $\Delta$, the limit
(where $F_{n}=F_{1}^{n}$ is the $n$-th iterate of $F_{1}$)
\[
h(z):=\lim\limits_{n\rightarrow\infty}\frac{1-F_{n}(z)}{1-F_{n}(0)}\
, \quad z\in\Delta,
\]
exists and is not constant (see \cite{E-S-V}). Moreover, for each
$t>0$, the function $h$ is the unique univalent solution of
Schr\"oder's functional equation
$$h(F_{t}(z))=\alpha^{t}h(z)$$ normalized by $h(0)=1$ (see
\cite{E-S-V} and \cite{E-G-R-S}). Hence,
$$h(G_{1}(z))=\lim\limits_{n\rightarrow\infty}\frac{1-F_{n}(G_{1}(z))}{1-F_{n}(0)}=\lim\limits_{n\rightarrow\infty}\frac{1-G_{1}(F_{n}(z))}{1-F_{n}(z)}\cdot\frac{1-F_{n}(z)}{1-F_{n}(0)}=\beta
h(z).$$ Therefore $$h(G_{1}(F_{t}(z)))=\beta
h(F_{t}(z))=\beta\alpha^{t}h(z)=\alpha^{t}h(G_{1}(z))=h(F_{t}(G_{1}(z)))
$$ for all $t\geq 0$ and $z\in\Delta$, and by the univalence of
$h$, $G_{1}$ commutes with $F_{t}$ for each $t\geq 0$.

Fix $t>0$, and denote by $\sigma$ the K{\oe}nigs function for
$S_{2}$: $$\sigma
(z):=\lim\limits_{n\rightarrow\infty}\frac{1-G_{n}(z)}{1-G_{n}(0)},
\, \, z\in\Delta.$$ Since the mapping $G_{1}$ is of hyperbolic
type, this limit exists and for each $s>0$, the function $\sigma$
is the unique univalent solution of Schr\"oder's functional
equation $$\sigma(G_{s}(z))=\beta^{s}\sigma(z)$$ normalized by
$\sigma(0)=1$. Hence,
$$\sigma(F_{t}(z))=\lim\limits_{n\rightarrow\infty}\frac{1-G_{n}(F_{t}(z))}{1-G_{n}(0)}=\lim\limits_{n\rightarrow\infty}\frac{1-F_{t}(G_{n}(z))}{1-G_{n}(z)}\cdot\frac{1-G_{n}(z)}{1-G_{n}(0)}=\alpha^{t}
\sigma(z).$$ Consequently,
$$\sigma(F_{t}(G_{s}(z)))=\alpha^{t}\sigma(G_{s}(z))=\alpha^{t}\beta^{s}\sigma(z)=\beta^{s}\sigma(F_{t}(z))=\sigma(G_{s}(F_{t}(z)))$$
for all $s>0$ and $z\in\Delta$, and by the univalence of $\sigma$
the semigroups commute.
\epr

\section{Semigroups of parabolic type}

For each $n=0,1,\ldots,$ we denote by $C_{A}^{n}(\tau)$,
$\tau\in\Delta$, the class of functions $F\in\Hol(\Delta
,\mathbb{C})$ which admit the representation
\begin{equation}\label{1zz}
F(z)=\sum_{k=0}^{n}a_{k}(z-\tau)^{k}+\gamma (z),
\end{equation}
where $\gamma\in\Hol(\Delta ,\mathbb{C})$ and
$\angle\lim\limits_{z\rightarrow\tau}\frac{\gamma
(z)}{(z-\tau)^{n}}=0$; and we say that $F\in C^{n}(\tau)$ when
this expansion holds as $z\rightarrow \tau$ unrestrictedly.

To proceed we need the following auxiliary result.

\begin{lemma}
Let $F,G\in \Hol(\Delta)$ be two commuting univalent parabolic
mappings and let $\tau =1$ be the Denjoy--Wolff point of $F$. If
one of the following conditions

(i) $F,G\in C^{2}(1)$, $F''(1)\neq 0$, $G''(1)\neq 0$;

(ii) $F,G\in C_{A}^{2}(1)$, $G''(1)\neq 0$, $\mathrm{Re}\,F''(1)>
0$;

(iii) $F,G\in C^{3}(1)$, $F''(1)=G''(1)=0$, $F'''(1)\neq 0$,
$G'''(1)\neq 0$

\noindent holds, then there exists a univalent function $\sigma\in
\Hol(\Delta,\mathbb{C})$ such that
\begin{equation}\label{1z}
\sigma\circ F=\sigma +1
\end{equation}
and
\begin{equation}\label{2z}
\sigma\circ G =\sigma +\lambda \quad  \mbox{with} \quad
\lambda\in\mathbb{C}, \, \, \, \lambda\neq 0.
\end{equation}
\end{lemma}

\pr Consider $z_{n}^{0}:=F_{n}(0)$ and
$\sigma_{n}(z):={\displaystyle\frac{F_{n}(z)-z_{n}^{0}}{z_{n+1}^{0}-z_{n}^{0}}},
\quad z\in\Delta.$ Then $\sigma_{n}\in\Hol(\Delta,\mathbb{C})$ and
the sequence $\{{\sigma_{n}}\}_{n=1}^{\infty}$ converges in the
compact-open topology to a certain holomorphic map
$\sigma\in\Hol(\Delta,\mathbb{C})$ such that (\ref{1z}) holds (by
Theorem 2.1 in \cite{C-M-P1}). Since $F$ is univalent in $\Delta$,
the solution $\sigma$ of Abel's equation (\ref{1z}) is also
univalent in $\Delta$.

Now we show that $\sigma$ satisfies (\ref{2z}). Denote $f=C\circ
F\circ C^{-1}$, $g=C\circ G\circ C^{-1}$,
$f,g\in\Hol(\mathbb{H},\mathbb{H})$,  where
$\mathbb{H}=\{z:\mathrm{Re}\,(z)>0\}$ and $C$ is the Cayley
transformation given by $C(z)=\frac{1+z}{1-z}$. Then $f$ and $g$
are commuting parabolic maps in $\Hol(\mathbb{H},\mathbb{H})$
having $\infty$ as their common Denjoy--Wolff point.

Denote $w_{0}:=C(0)=1$, $w_{n}^{0}:=f_{n}(1)=C(z_{n}^{0})$,
$$w_{n}:=f_{n}(w),\quad w_{n}\in\mathbb{H},$$ and $$h_{n}(w):=\frac{w_{n}-w_{n}^{0}}{w_{n+1}^{0}-w_{n}^{0}}\,,\quad
\, \, w\in\mathbb{H}.$$ Then $h_{n}\in\Hol(\mathbb{H},\mathbb{C})$
and the sequence $\{{h_{n}}\}_{n=1}^{\infty}$ converges in the
compact open topology to a holomorphic function
$h\in\Hol(\mathbb{H},\mathbb{C})$ such that $h\circ f=h+1$ and
$\sigma=h\circ C$ (see \cite{C-M-P1}).

Suppose that (i) holds. Then the following expansions of $f$ and
$g$ at $\infty$ are satisfied (see \cite{B-S}):
\begin{equation}\label{3z}
f(w)=w+F''(1)+\gamma_{f}(w), \, \,
\lim\limits_{w\rightarrow\infty}\gamma_{f}(w)=0
\end{equation}
\noindent and
\begin{equation}\label{4z}
g(w)=w+G''(1)+\gamma_{g}(w), \, \,
\lim\limits_{w\rightarrow\infty}\gamma_{g}(w)=0.
\end{equation}
Hence,
$$h_{n}(g(w))=\frac{f_{n}(g(w))-w_{n}^{0}}{w_{n+1}^{0}-w_{n}^{0}}=\frac{g(f_{n}(w))-w_{n}^{0}}{w_{n+1}^{0}-w_{n}^{0}}$$
$$=\frac{w_{n}+G''(1)+\gamma_{g}(w_{n})-w_{n}^{0}}{w_{n+1}^{0}-w_{n}^{0}}=\frac{w_{n}-w_{n}^{0}}{w_{n+1}^{0}-w_{n}^{0}}+\frac{G''(1)+\gamma_{g}(w_{n})}{w_{n+1}^{0}-w_{n}^{0}}$$
$$=h_{n}(w)+\frac{G''(1)+\gamma_{g}(w_{n})}{F''(1)+\gamma_{f}(w_{n})}\cdot\frac{F''(1)+\gamma_{f}(w_{n})}{w_{n+1}^{0}-w_{n}^{0}}\,
.$$
Letting $n\rightarrow\infty$, we obtain
\begin{equation}\label{4d}
h(g(w))-h(w)=\frac{G''(1)}{F''(1)}\cdot \lim\limits_{n\rightarrow
\infty}\frac{F''(1)+\gamma_{f}(w_{n})}{w_{n+1}^{0}-w_{n}^{0}}\,.
\end{equation}
Repeating this calculation with $f$ instead of $g$, we find that
$$h(f(w))=h(w)+\lim\limits_{n\rightarrow
\infty}\frac{F''(1)+\gamma_{f}(w_{n})}{w_{n+1}^{0}-w_{n}^{0}}\,.$$
At the same time, $h\circ f=h+1$. Hence
${\displaystyle\lim\limits_{n\rightarrow
\infty}\frac{F''(1)+\gamma_{f}(w_{n})}{w_{n+1}^{0}-w_{n}^{0}}=1}$.

Rewrite (\ref{4d}) as follows: $$h(g(w))-h(w)=\lambda, \quad
\mbox{where} \quad \lambda=\frac{G''(1)}{F''(1)}\neq 0 \quad
\mbox{and} \quad w\in\mathbb{H}.$$ Substituting $h=\sigma\circ
C^{-1}$ and $g=C\circ G \circ C^{-1}$ in the last equality we get
(\ref{2z}).

If (ii) holds, then Theorem 14 in \cite{C-M-P} implies that for
each $z\in\Delta$, the sequence $\{F_{n}(z)\}_{n=1}^{\infty}$
converges to 1 (and, consequently, $\{w_{n}\}$ converges to
$\infty$) nontangentially. So, in this case, one can repeat the
proof of item (i), replacing the unrestricted limits in (\ref{3z})
and (\ref{4z}) by the angular limits.

Suppose now that (iii) holds. Then the following expansions of $f$
and $g$ at $\infty$ hold (see \cite{B-S}):
\begin{equation}\label{5z}
f(w)=w-\frac{2}{3}\frac{F'''(1)}{w+1}+\Gamma_{f}(w), \, \,
\lim\limits_{w\rightarrow\infty}\Gamma_{f}(w)w=0
\end{equation}
\noindent and
\begin{equation}\label{6z}
g(w)=w-\frac{2}{3}\frac{G'''(1)}{w+1}+\Gamma_{g}(w), \, \,
\lim\limits_{w\rightarrow\infty}\Gamma_{g}(w)w=0.
\end{equation}
Therefore
$$h_{n}(f(w))=\frac{f(w_{n})-w_{n}^{0}}{w_{n+1}^{0}-w_{n}^{0}}=\frac{w_{n}-w_{n}^{0}}{w_{n+1}^{0}-w_{n}^{0}}+\frac{-\frac{2}{3}\frac{F'''(1)}{w_{n}+1}+\Gamma_{f}(w_{n})
 }{w_{n+1}^{0}-w_{n}^{0}}$$
$$=h_{n}(w)+\frac{-\frac{2}{3}\frac{F'''(1)}{w_{n}+1}+\Gamma_{f}(w_{n})
}{w_{n+1}^{0}-w_{n}^{0}}\,.$$ Letting $n\rightarrow\infty$, we
obtain $$h(f(w))=h(w)+\lim\limits_{n\rightarrow
\infty}\frac{-\frac{2}{3}\frac{F'''(1)}{w_{n}+1}+\Gamma_{f}(w_{n})
}{w_{n+1}^{0}-w_{n}^{0}}\,.$$ On the other hand, $h(f(w))=h(w)+1.$
Hence,
\begin{equation}\label{7z}
\lim\limits_{n\rightarrow
\infty}\frac{-\frac{2}{3}\frac{F'''(1)}{w_{n}+1}+\Gamma_{f}(w_{n})}{w_{n+1}^{0}-w_{n}^{0}}=1.
\end{equation}
Now using (\ref{6z}), we find
$$h_{n}(g(w))=\frac{g(w_{n})-w_{n}^{0}}{w_{n+1}^{0}-w_{n}^{0}}=\frac{w_{n}-\frac{2}{3}\frac{G'''(1)}{w_{n}+1}+\Gamma_{g}(w_{n})-w_{n}^{0}
}{w_{n+1}^{0}-w_{n}^{0}}$$
$$=h_{n}(w)+\frac{-\frac{2}{3}G'''(1)+\Gamma_{g}(w_{n})(w_{n}+1)}{-\frac{2}{3}F'''(1)+\Gamma_{f}(w_{n})(w_{n}+1)
}\cdot\frac{-\frac{2}{3}\frac{F'''(1)}{w_{n}+1}+\Gamma_{f}(w_{n})}{w_{n+1}^{0}-w_{n}^{0}}\,.$$
Letting $n\rightarrow\infty$ and using (\ref{7z}), we get
$$h(g(w))-h(w)=\lambda, \, \, w\in\mathbb{H}, \, \, \mbox{where}
\, \, \lambda=\frac{G'''(1)}{F'''(1)}\neq 0.$$ Consequently,
$\sigma\circ G-\sigma =\lambda$. \epr

Following \cite{C-M-P1}, we say that the function $\sigma$
mentioned in the lemma is {\it the K{\oe}nigs intertwining
function} associated with $F$ with respect to $z_{0}=0$.

\begin{rem}
The function $\sigma$ in Lemma 2 is completely determined by the
function $F$. It does not depend on $G$. So if the conditions of
the lemma hold for the same function $F$ and another function
$G_{1}\in\Hol(\Delta)$, then we have the equality
$$\sigma\circ G_{1}=\sigma+\lambda_{1}$$ with the same function
$\sigma$ and a constant $\lambda_{1}\neq 0$.
\end{rem}

\begin{pred0}[parabolic case]
Let $S_{1}=\left\{F_{t}\right\}_{t\geq 0}$ and
$S_{2}=\left\{G_{t}\right\}_{t\geq 0}$ be two continuous
semigroups on $\Delta$ generated by $f$ and $g$, respectively, and
let $F_{1}\circ G_{1}=G_{1}\circ F_{1}$.

Suppose that $\tau=1$ is the boundary null point of $f$ such that
$f'(1)=0$. If $S_{1},S_{2}\subset C^0(1)$ and one of the following
conditions

(i) $f,g\in C^{2}(1)$, $f''(1)\neq 0$, $g''(1)\neq 0$;

(ii) $f,g\in C^{3}(1)$, $f''(1)=g''(1)=0$

\noindent holds, then the semigroups commute.
\end{pred0}

\pr
Since $\tau$ is a boundary null point of $f$ and $f'(\tau)=0$, it
is the common Denjoy--Wolff point of the semigroup $S_{1}$. The
commutativity of $F_{1}$ and $G_{1}$ implies that $\tau$ is the
Denjoy--Wolff point of $G_{1}$ (see \cite{Beh}) and, consequently,
$\tau$ is also the common Denjoy--Wolff point of the semigroup
$S_{2}$ .

If (ii) holds and, in addition, either $f'''(1)=0$ or $g'''(1)=0$,
then by Corollary 1 we have that either $F_{t}\equiv I$ or
$G_{t}\equiv I$, respectively, and therefore the semigroups
commute. Suppose that $f'''(1)\neq 0$ and $g'''(1)\neq 0$ in (ii).
Then by Remark 3 above, one can replace conditions (i) and (ii) by

(i') $F_{t}, G_{t}\in C^{2}(1)$, $F''_{t}(1)\neq 0$,
$G''_{t}(1)\neq 0$ for all $t>0$;

(ii') $F_{t}, G_{t}\in C^{3}(1)$, $F''_{t}(1)=G''_{t}(1)=0$,
$F'''_{t}(1)\neq 0$, $G'''_{t}(1)\neq 0$, $t>0$.

By our assumption, $F_{1}\circ G_{1}=G_{1}\circ F_{1}$. Moreover,
$F_{1}\circ F_{t}=F_{t}\circ F_{1}$ for all $t\geq 0$. Therefore
Lemma 2 implies that there exists the K{\oe}nigs intertwining map
$\sigma$ for $F_{1}$ with respect to $z_{0}=0$, which satisfies
\begin{equation}\label{8z}
\sigma(F_{1}(z))=\sigma(z)+1, \, \, z\in\Delta,
\end{equation}
\begin{equation}\label{9z}
\sigma(G_{1}(z))-\sigma(z)=\lambda,  \, z\in\Delta, \, \, \,
\mbox{for some} \,  \lambda\neq 0,
\end{equation}
\noindent and
\begin{equation}\label{10z}
\sigma(F_{t}(z))-\sigma(z)=\beta (t), \, \, t>0, \, \, z\in\Delta,
\end{equation}
where $\beta (t)\neq 0$ for all $t>0$.

Furthermore, $F_{1}\circ G_{1}=G_{1}\circ F_{1}$ and $G_{1}\circ
G_{s}=G_{s}\circ G_{1}$ for all $s\geq 0$. Hence, by Lemma 2,
there exists the K{\oe}nigs intertwining map $\widetilde{\sigma}$
for $G_{1}$ with respect to $z_{0}=0$, which satisfies
\begin{equation}\label{11z}
\widetilde{\sigma}(G_{1}(z))=\widetilde{\sigma}(z)+1, \, \,
z\in\Delta,
\end{equation}
\begin{equation}\label{12z}
\widetilde{\sigma}(F_{1}(z))-\widetilde{\sigma}(z)=\widetilde{\lambda},
 \, z\in\Delta, \, \mbox{for some} \,  \widetilde{\lambda}\neq 0,
\end{equation}
\noindent and
\begin{equation}\label{13z}
\widetilde{\sigma}(G_{s}(z))-\widetilde{\sigma}(z)=\widetilde{\beta}
(s), \, \, s>0, \, \, z\in\Delta,
\end{equation}
where $\widetilde{\beta}(s)\neq 0$ for all $s>0$.

Assume that at least one of the mappings $F_{1}, \, G_{1}$ (for
example, $G_{1}$) is of nonautomorphic type. (Note that if (ii')
holds then for each $t>0$, $G_{t}$ and $F_{t}$ are of
nonautomorphic type by Theorem 4.4 in \cite{Shap}.) It follows
from (\ref{12z}) and (\ref{13z}) that
\begin{equation}\label{14z}
\frac{\widetilde{\beta}(s)}{\widetilde{\lambda}}(\widetilde{\sigma}(F_{1}(z))-\widetilde{\sigma}(z))=\widetilde{\sigma}(G_{s}(z))-\widetilde{\sigma}(z).
\end{equation}
Rewrite (\ref{9z}) in the form
\begin{equation}\label{15z}
\frac{1}{\lambda}\sigma(G_{1}(z))=\frac{1}{\lambda}\sigma(z)+1.
\end{equation}
By Theorem 3.1 in \cite{C-M-P1}, equalities (\ref{11z}) and
(\ref{15z}) imply that
$\frac{1}{\lambda}\sigma=\widetilde{\sigma}+const.$, and so
(\ref{14z}) is equivalent to
\begin{equation}\label{16z}
\frac{\widetilde{\beta}(s)}{\widetilde{\lambda}}(\sigma(F_{1}(z))-\sigma(z))=\sigma(G_{s}(z))-\sigma(z)
\end{equation}
or, by (\ref{8z}),
\begin{equation}\label{17z}
\frac{\widetilde{\beta}(s)}{\widetilde{\lambda}}=\sigma(G_{s}(z))-\sigma(z).
\end{equation}
Since the right-hand sides in (\ref{10z}) and (\ref{17z}) are
differentiable in $t$ and $s$, respectively, $\beta$ and
$\widetilde{\beta}$ are differentiable too. Hence,
$$\beta'(t)=\sigma '(F_{t}(z))\cdot \frac{\partial
F_{t}(z)}{\partial t} \, \, \, \,  \mbox{and} \, \, \, \,
\frac{\widetilde{\beta}'(s)}{\widetilde{\lambda}}=\sigma
'(G_{s}(z))\cdot \frac{\partial G_{s}(z)}{\partial s}.$$

Letting $t\rightarrow 0^{+}$ and $s\rightarrow 0^{+}$ in these
equalities, we obtain $$\beta'(0)=-\sigma '(z)\cdot f(z) \, \, \,
\,  \mbox{and} \, \, \, \,
\frac{\widetilde{\beta}'(0)}{\widetilde{\lambda}}=-\sigma
'(z)\cdot g(z),$$ where $f$ and $g$ are generators of the
semigroups $\left\{F_{t}\right\}_{t\geq 0}$ and
$\left\{G_{t}\right\}_{t\geq 0}$, respectively.

Since $\sigma$ is univalent on $\Delta$, the derivative $\sigma
'(z)\neq 0$ for all $z\in\Delta$. Moreover, because the common
Denjoy--Wolff point of $S_{1}$ and $S_{2}$ belongs to the boundary
$\partial\Delta$, the generators $f$ and $g$ do not vanish on
$\Delta$. Therefore $$f(z)=ag(z), \, \, \mbox{where} \, \,
a=\frac{\widetilde{\lambda}\beta '(0)}{\widetilde{\beta} '(0)},$$
and by \cite{E-L-S-T}, the semigroups commute.

Now let the mappings $F_{1}$ and $G_{1}$ be both of automorphic
type. Note that in this case $F''_{1}(1)$ and $G''_{1}(1)$ cannot
be zero and so we assume that (i') holds.

We have already seen in the proof of Lemma 2 that
\begin{equation}\label{18z}
\sigma(G_{1}(z))-\sigma (z)=\frac{G''_{1}(1)}{F''_{1}(1)}.
\end{equation}

Since $\mathrm{Re}\, F''_{1}(1)=0$ and $\mathrm{Re}\,
G''_{1}(1)=0$ (see Theorem 4.4 in \cite{B-S}), it follows that
$\frac{G''_{1}(1)}{F''_{1}(1)}\in\mathbb{R}\setminus\{0\}$.
Moreover, by Theorem 1, $$F''_{t}(1)=-\alpha t \quad \mbox{and}
\quad G''_{t}(1)=-\widetilde{\alpha} t, \, \, t>0, $$ where
$\alpha=f''(1)\neq 0$ and $\widetilde{\alpha}=g''(1)\neq 0$. So
equality (\ref{18z}) has the form
\begin{equation}\label{19z}
\sigma (G_{1}(z))-\sigma (z)=p, \quad\mbox{where} \quad
p:=\frac{\widetilde{\alpha}}{\alpha}.
\end{equation}
On the other hand,
\begin{equation}\label{20z}
\sigma (F_{t}(z))-\sigma
(z)=\frac{F''_{t}(1)}{F''_{1}(1)}=\frac{\alpha t}{\alpha}=t \quad
 \mbox{for all} \quad t\geq 0.
\end{equation}

First we suppose that $p>0$. From (\ref{19z}) and (\ref{20z}) we
have $\sigma(G_{1}(z))=\sigma(F_{p}(z)), \, \, z\in\Delta, $ and
by the univalence of $\sigma$ on $\Delta$, $G_{1}(z)=F_{p}(z)$ for
all $z\in\Delta$. Hence, $G_{1}\circ F_{t}=F_{t}\circ G_{1}$ for
all $t\geq 0$.

Fix $t>0$ and repeat these considerations with $G_{1}$, $F_{t}$,
$G_{s}$ and $\widetilde{\sigma}$ instead of $F_{1}$, $G_{1}$,
$F_{t}$ and $\sigma$, respectively. Namely,
$$\widetilde{\sigma}(F_{t}(z))-\widetilde{\sigma}(z)=
\frac{F''_{t}(1)}{G''_{1}(1)}=\frac{\alpha t
 }{\widetilde{\alpha}}>0$$
\noindent and
$$\widetilde{\sigma}(G_{s}(z))-\widetilde{\sigma}(z)=
\frac{G''_{s}(1)}{G''_{1}(1)}=s \quad \mbox{for all} \quad s>0. $$

Denote $\widetilde{s}:=\frac{\alpha t}{\widetilde{\alpha} }>0$.
Then $ \widetilde{\sigma} (F_{t}(z))=\widetilde{\sigma}
(G_{\widetilde{s}}(z)), \, \, z\in\Delta$. By the univalence of
$\widetilde{\sigma}$ on $\Delta$ we have
$F_{t}(z)=G_{\widetilde{s}}(z)$. Therefore $G_{s}\circ
F_{t}=F_{t}\circ G_{s}$ for all $s>0$. Since $t>0$ is arbitrary,
it follows that the semigroups $S_{1}=\left\{F_{t}\right\}_{t\geq
0}$ and $S_{2}=\left\{G_{s}\right\}_{s\geq 0}$ commute.

Let now $p<0$. Then by (\ref{20z}), $\sigma (F_{-p}(z))-\sigma
(z)=-p$ for all $z\in\Delta$. Hence, by (\ref{19z}),
\[
\sigma (F_{-p}(G_{1}(z)))-\sigma (G_{1}(z))=\sigma (z)-\sigma
(G_{1}(z)),\quad z\in\Delta,
\]
and, therefore,
\[
\sigma (F_{-p}(G_{1}(z)))=\sigma (z), \quad z\in\Delta.
\]
By the univalence of $\sigma$ on $\Delta$, $F_{-p}(G_{1}(z))=z$.
Consequently, $F_{-p}=G_{1}^{-1}$ on $G_{1}(\Delta)$. Since
$F_{-p}\in\Hol(\Delta)$, $G_{1}^{-1}$ is well defined on $\Delta$
and so $G_{1}$, as well as $F_{-p}$, are an automorphisms of
$\Delta$. Therefore, by Proposition 5,
$\left\{F_{t}\right\}_{t\geq 0}$ is a semigroup of automorphisms.
Consequently, it can be extended to a group
$S_{F}=\left\{F_{t}\right\}_{t\in\mathbb{R}}$ and
$G_{1}=F_{p}^{-1}=F_{-p}\in S_{F}$. In particular, $G_{1}\circ
F_{t}=F_{t}\circ G_{1}$ for all $t\geq 0$.

Fix $t>0$. In a similar way, using the commutativity of $F_{t}$
and $G_{1}$, one can show that the semigroup
$\left\{G_{s}\right\}_{s\geq 0}$ can be extended to a group
$S_{G}=\left\{G_{s}\right\}_{s\in\mathbb{R}}$ and that $F_{t}\circ
G_{s}=G_{s}\circ F_{t} $ for all $s,t\in\mathbb{R}$. \epr

\begin{rem}
{\rm Note in passing that the proof of Theorem 4 implies the
following interesting fact:}

Let $S_{1}=\left\{F_{t}\right\}_{t\geq 0}$ be a continuous
semigroup of parabolic type on $\Delta$ generated by $f$ with the
Denjoy--Wolff point $\tau=1$, and let $G$ be a holomorphic
self-mapping of\  $\Delta$ such that $F_{1}\circ G=G\circ F_{1}$.
If $\, f,G\in C^{2}(1)$ and $S_{1}\subset C^{0}(1)$, then the
condition $f''(1)\cdot G''(1)>0$ implies that $S_1$ can be
extended to a group of parabolic automorphisms of $\Delta$ and
$G\in S_1$, hence $G$ commute with all elements $F_{t},\ t>0$.
\end{rem}

\begin{rem}
{\rm Note also that if in the assumptions of Theorem~4,}
$S_{1}=\left\{F_{t}\right\}_{t\geq 0}$ and
$S_{2}=\left\{G_{t}\right\}_{t\geq 0}$ are both groups of
parabolic automorphisms of $\Delta$, then condition (i) of the
theorem holds automatically, so the commutativity of $F_{1}$ and
$G_{1}$ implies that $S_{1}$ and $S_{2}$ commute.
\end{rem}
%
%

\vspace{5mm}

{\bf Acknowledgment.} The third author was partially supported by
the Fund for the Promotion of Research at the Technion and by the
Technion VPR Fund - B. and G. Greenberg Research Fund (Ottawa).

\end{document}